\documentstyle[amssymb,12pt]{article}
\newcommand{\qed}{{\hspace{\fill} \lower .5mm \hbox{$\Box$}}}

\begin{document}

\begin{center}
{TENSOR PRODUCTS OF MAXIMAL ABELIAN SUBALGBERAS OF C*-ALGEBRAS.}
\footnote{Version of November 23, 2007}

\bigskip
SIMON WASSERMANN

{\footnotesize\it Department of Mathematics, University of
Glasgow, Glasgow G12 8QW, United Kingdom}
\end{center}

\bigskip
{\bf Abstract.} It is shown that if $C_1$ and $C_2$ are maximal
abelian self-adjoint subalgebras (masas) of C*-algebras $A_1$ and
$A_2$, respectively, then the completion $C_1\otimes C_2$ of the
algebraic tensor product $C_1\odot C_2$ of $C_1$ and $C_2$ in any
C*-tensor product $A_1\otimes_{\beta} A_2$ is maximal abelian
provided that $C_1$ has the extension property of Kadison and
Singer and $C_2$ contains an approximate identity for $A_2$.
Examples are given to show that this result can fail if the
conditions on the two masas do not both hold. This gives an answer
to a long-standing question, but leaves open some other
interesting problems, one of which turns out to have a potentially
intriguing implication for the Kadison-Singer extension problem.


\bigskip {\bf 1. Introduction.} If $A_1$ and $A_2$ are C*-algebras
with centres $Z_1$ and $Z_2$ respectively, it was shown by Richard
Haydon and the author [{\bf 6}] that the centre of the minimal
C*-tensor product $A_1\otimes_{min}A_2$ is just the closure of the
algebraic tensor product $Z_1\odot Z_2$ in $A_1\otimes_{min}A_2$.
This closure is naturally isomorphic to $Z_1\otimes_{min} Z_2$.
The result was also shown to follow from a more general slice map
result [{\bf 12}, Theorem 4] and the analogous result for
arbitrary C*-tensor products was subsequently established by
Archbold using the Dixmier approximation property [{\bf 3}]. Batty
later gave a neat alternative proof of Archbold's general result
[{\bf 5}].

Analogous questions arise for maximal abelian self-adjoint
C*-subalgebras (or {\it masas}) of C*-algebras. If $C_1$ and $C_2$
are masas of $A_1$ and $A_2$, respectively, then the closure of
the algebraic tensor product $C_1\odot C_2$ in any C*-completion
$A_1\otimes_{\beta}A_2$ of $A_1\odot A_2$ is naturally isomorphic
to $C_1\otimes_{min}C_2$ since the algebras in the tensor product
are abelian. The slice map result [{\bf 12}, Theorem 4] implies
that $C_1\otimes_{min}C_2$ is again a masa in $A_1\otimes_{min}
A_2$, and it is natural to ask whether $C_1\otimes_{min}C_2$ is a
masa in $A_1\otimes_{\beta}A_2$ for any C*-norm $\beta$ on
$A_1\odot A_2$ other than the minimal norm when $A_1\odot A_2$ has
more than one C*-norm.

Although this question was originally raised in [12], up to now
little progress seems to have been made. In this paper we give a
solution to the problem. There are two main results. The first, a
positive one, states that if one of the masas $C_1$, $C_2$
possesses the Kadison-Singer extension property and the other
contains an approximate identity for the ambient algebra, then the
question has a positive answer. The second, a negative answer to
the general question, is a pair of example of masas $C_1$ and
$C_2$ of C*-algebras $A_1$ and $A_2$, respectively, such that
$C_1\otimes_{min} C_2$ is not maximal abelian in $A_1\otimes_{max}
A_2$ in each case. Connections with the Kadison-Singer extension
problem [{\bf 9}] of whether $\ell^\infty({\Bbb N})$ has the
extension property relative to $B(\ell^2({\Bbb N})$ will be
discussed in the final section.

\bigskip {\bf 2. Masas with the extension property.} A masa $C$
of a C*-algebra $A$ is said to have the {\it extension property}
(see [{\bf 9}], [{\bf 2}], [{\bf 4}]) if

(i) any pure state (i.e.\ character) of $C$ has a unique pure
state extension to $A$ and (when $A$ is non-unital)

(ii) no pure state of $A$ annihilates $C$.

\noindent Condition (ii) is well-known to be equivalent to the
condition (see [{\bf 1}, Proof of Lemma 2.32])

(ii)$'$ $C$ contains an increasing approximate identity for $A$.

\noindent It is a straightforward consequence of the Krein-Milman
theorem that any pure state of a masa $C$ with the extension
property has a unique state extension to $A$. An alternative
characterisation [{\bf 4}] of the extension property for unital
$A$ states that, if ${\cal U}(C)$ denotes the unitary group of
$C$, then for any $x\in A$ the intersection $C\cap
\overline{co}\{uxu^*:u\in {\cal U}(C)\}$ of $C$ with the closed
convex hull $\overline{co}\{uxu^*:u\in {\cal U}(C)\}$ contains
exactly one point.

\bigskip {\sc Examples.} 1. In the reduced C*-algebra
$C^*_r({\Bbb F}_2)$ of the free group on two generators, with $u$
and $v$ the canonical unitary generators, the abelian
C*-subalgebras generated by $u$ and $v$, respectively, are masas
with the extension property [{\bf 4}, Example (i)]. Moreover
$C^*_r({\Bbb F}_2)$ is not nuclear [{\bf 11}] (see also [{\bf
14]}).

\bigskip 2. No non-atomic masa of $B(\ell^2({\Bbb N}))$
has the extension property [{\bf 9}].

\bigskip In what follows the minimal C*-tensor product of two C*-algebras
will be denoted by $A\otimes B$ when at least one of $A$ and $B$
is abelian. The following well-known factorisation result for
states on a tensor product will be required in the proof of
Theorem 2.

\bigskip{\sc Lemma 1} [{\bf 5}, Lemma 3]. {\it Let
$\varphi$ be a state on $A_1\otimes_{\beta}A_2$ such that the
restriction $\varphi_1$ of $\varphi$ to $A_1$ is a pure state of
$A_1$. Then $\varphi=\varphi_1\otimes \varphi_2$ for some state
$\varphi_2$ on $A_2$.}

\bigskip In the proof of the following result, which is analogous to that of
[{\bf 5}, Theorem 4], the unitization $\tilde{A}$ of $A$ will be
taken to be the subalgebra $A + {\Bbb C}.1$ of the multiplier
algebra $M(A)$, so that $\tilde{A} = A$ if $A$ is itself unital.
If $\beta$ is a C*-norm on $A_1\odot A_2$, $\tilde{\beta}$ will
denote the unique C*-norm on $\tilde{A}_1\odot \tilde{A}_2$
extending $\beta$ (see [{\bf 8}]).

\bigskip {\sc Theorem 2.} {\it Let $A_1$ and $A_2$ be
C*-algebras with $A_2$ unital and let $C$ be a masa of $A_1$ with
the extension property. Then in any C*-completion
$A_1\otimes_{\beta} A_2$,
   \[(C\otimes1)^c = C\otimes A_2,\]
where $(C\otimes 1)^c = \{x\in A_1\otimes_{\beta} A_2: x(c\otimes
1) = (c\otimes 1)x \mbox{ for all }c\in C\}$. }

\bigskip
{\it Proof.} 1. Assume first that $A_1$ is unital. If
$\Phi:A_1\otimes_{\beta}A_2 \to A_1\otimes_{min}A_2$ is the
canonical homomorphism, its restriction to $C\otimes A_2$ is an
isomorphism and $C\otimes A_2$ can be identified with its image in
$A_1\otimes_{min}A_2$ under $\Phi$. With this identification
$\Phi|_{C_1\otimes A_2}$ is just the identity map. If $x\in
(C\otimes 1)^c$, then $[x,C\otimes 1] = 0$, which implies that
$[\Phi(x),C\otimes 1] = 0$. By [{\bf 12}, Theorem 4], $\Phi(x)\in
C\otimes A_2$. Replacing $x$ by $(x-\Phi(x))^*(x-\Phi(x))$, it is
sufficient to show that if $x\geq 0$ and $\Phi(x) = 0$ then $x=0$.

With these assumptions $\|x\|\in $ Sp$(x)$, $\|x\|1 - x$ is
singular and the closed left ideal $I$ of $D=(C\otimes 1)^c$
generated by $\|x\|1-x$ is proper. Let
   \[J = \{c\in C: c\otimes 1 \in I\}.\]
Since $I(C\otimes 1) \subseteq I$, $J$ is a proper closed
two-sided ideal of $C$ and since $J\otimes 1 = (C\otimes 1)\cap
I$, there are canonical isometric isomorphisms
   \[C/J \cong (C\otimes 1)/(J\otimes 1) \cong
(C\otimes 1 + I)/I,\] by [{\bf 10}, 1.17.6]. Let $\chi$ be a
character of $C$ which annihilates $J$. Via these isomorphisms,
$\chi$ corresponds to a linear functional $\varphi$ on $(C\otimes
1 + I)/I$ such that $\|\varphi\| = 1$ and $\varphi(1+I) = 1$. By
the Hahn-Banach theorem $\varphi$ extends to a linear functional
of norm $1$ on $D/I$ which, when composed with the quotient map,
gives a state $\bar{\varphi}$ on $D$ such that
$\bar{\varphi}(c\otimes 1) = \chi(c)$ for $c\in C$. Let $\psi$ be
an extension of $\bar{\varphi}$ to a state on
$A_1\otimes_{\beta}A_2$. Letting $\psi_1$ be the restriction of
$\psi$ to $A_1$, so that $\psi_1(a) = \psi(a\otimes 1)$, $\psi_1$
is a pure state since $C$ has the extension property. By Lemma 1,
$\psi=\psi_1\otimes \psi_2$ for some state $\psi_2$ on $A_2$. Now
$\psi(\|x\|1 - x) = 0$, since $\psi|_I = 0$, and $\psi(x) =
(\psi_1\otimes \psi_2)(\Phi(x)) = 0$. Thus $\|x\| = \psi(\|x\|1 -
x) = 0$, which implies that $x = 0$, as required.

2. If $A_1$ is not unital, let $\tilde{A}_1$ be the unitisation of
$A_1$ and let $\tilde{C} = C + {\Bbb C}1$. Then $\tilde{C}$ is a
masa in $\tilde{A}_1$. To see that $\tilde{C}$ has the extension
property in $\tilde{A}_1$, let $f$ be a pure state of $\tilde{C}$
and let $\bar{f}$ be a pure state extension of $f$ to
$\tilde{A}_1$. If $f$ is the unique pure state annihilating $C$,
then $g=\bar{f}|_{A_1}=0$, since otherwise $g$ would be a pure
state of $A_1$ which annihilated $C$. In this case $\bar{f}$ is
the unique pure state of $\tilde{A}_1$ which annihilates $A_1$. If
$f|_C\neq 0$, then $f|_C$ is a pure state of $C$ and
$\bar{f}|_{A_1}$ is a pure state extension of $f|_C$, hence
uniquely determined by $f$. Since $\bar{f}$ is uniquely determined
by its restriction to $A_1$, it follows that $\bar{f}$ is uniquely
determined by $f$.

For $x\in A_1\otimes_{\beta}A_2$, if $x\in (C\otimes 1)^c$ then
$x\in (\tilde{C}\otimes
1)^c\subseteq\tilde{A}_1\otimes_{\tilde{\beta}}A_2$, which implies
by part 1 that $x\in \tilde{C}\otimes A_2$. Let $\varphi$ be the
state on $\tilde{A}_1$ which annihilates $A_1$. Then $\chi =
\varphi|_{\tilde{C}}$ is the character of $\tilde{C}$ which
annihilates $C$ and the kernel of the map $\chi\otimes id_{A_2}$
is just $C\otimes A_2$. Since $0=(\varphi\otimes id_{A_2})(x) =
(\chi\otimes id_{A_2})(x)$, it follows that $x\in C\otimes A_2$ as
required. \qed

\bigskip {\sc Corollary 3.} {\it Let $A_1$ and $A_2$ be
C*-algebras and let $C_1$ and $C_2$ be masas of $A_1$ and $A_2$,
respectively, such that $C_1$ has the extension property and $C_2$
contains an approximate identity for $A_2$ if $A_2$ is not unital.
Then $C_1\otimes C_2$ is a masa of $A_1\otimes_{\beta}A_2$ for any
C*-norm $\beta$ on $A_1\odot A_2$. Moreover $C_1\otimes C_2$ has
the extension property if and only if $C_2$ does.}

\bigskip
{\it Proof.} Assume that $C_1$ has the extension property. If
$A_2$ is unital, it is immediate from Theorem 2 that $(C_1\otimes
C_2)^c \subseteq (C_1\otimes 1)^c = C_1\otimes A_2$. If $A_2$ is
non-unital and $\{e_{\lambda}\}$ is an approximate identity of
$A_2$ in $C_2$, for $x\in (C_1\otimes C_2)^c$, $[x, c\otimes
e_{\lambda}] = 0$ for any $c\in C_1$ and any $\lambda$. Since
$x(c\otimes 1) = \lim_{\lambda}x(c\otimes e_{\lambda})$ and
$(c\otimes 1)x = \lim_{\lambda}(c\otimes e_{\lambda})x$, it
follows that, as an element of
$A_1\otimes_{\tilde{\beta}}\tilde{A}_2$, $x$ lies in $(C_1\otimes
1)^c$, which equals $C_1\otimes_{\tilde{\beta}} \tilde{A_2}$, the
closure of $C_1\odot \tilde{A}_2$ in $A_1\otimes_{\tilde{\beta}}
\tilde{A}_2$, by Theorem 2. Since $C_1$ is abelian,
$C_1\otimes_{\tilde{\beta}} \tilde{A_2}$ identifies naturally with
$C_1\otimes \tilde{A_2}$ and so $x\in (C_1\otimes \tilde{A_2})\cap
(A_1\otimes_{\beta} A_2)$. Let $\Psi$ denote the canonical
homomorphism from $A_1\otimes_{\tilde{\beta}}\tilde{A}_2$ to
$A_1\otimes_{\min}\tilde{A}_2$. Then $\Psi|_{C_1\otimes
\tilde{A}_2}$ is the identity map on $C_1\otimes \tilde{A}_2$,
since $C_1\otimes_{\tilde{\beta}} \tilde{A}_2$ and $C_1\otimes
\tilde{A}_2$ are naturally identified, and $\Psi(x)$ is in
$(C_1\otimes \tilde{A}_2)\cap(A_1\otimes_{min}A_2)$, which equals
$C_1\otimes A_2$ by [{\bf 12}, Cor.\ 5] (or a simple slice map
argument), since $C_1$ is abelian. Thus $x\in
C_1\otimes_{\beta}A_2$. If $X$ is the spectrum of $C_1$, there is
a natural isomorphism $C_1\otimes A_2\cong C_0(X,A_2)$. Via this
isomorphism $C_1\otimes C_2$ identifies with $C_0(X,C_2)$ and $x$
identifies with a function $\overline{x}$ in $C_0(X, A_2)$. For
$\chi\in X$ and $c\in C_2$, there is a function $f\in C_0(X, C_2)$
such that $f(\chi) = c$. Since $x\in (C_1\otimes C_2)^c$,
$[\overline{x}(\chi), c] = [\overline{x}, f](\chi) =  0$, which
implies that $\overline{x}(\chi)\in C_2$ for $\chi\in X$ and so
$\overline{x}\in C_0(X,C_2)$. Thus $x\in C_1\otimes C_2$, as
required.

Now assume that $C_1$ and $C_2$ have the extension property and
let $\varphi$ be a character of $C_1\otimes C_2$. Then
$\varphi=\varphi_1\otimes\varphi_2$ for suitable characters
$\varphi_1$ and $\varphi_2$ of $C_1$ and $C_2$, respectively. Let
$\psi_i$ be the unique pure state extension of $\varphi_i$ to
$A_i$ for $i=1,2$ and let $\psi$ be the pure state $\psi_1\otimes
\psi_2$ of $A_1\otimes_{\beta}A_2$. If $\bar{\varphi}$ is any pure
state extension of $\varphi$ to $A\otimes_{\beta}A_2$, and
$\bar{\varphi}_1$ and $\bar{\varphi}_2$ are the restrictions of
$\bar{\varphi}$ to $A_1$ and $A_2$, respectively, then
$\bar{\varphi}_i|_{C_i}=\varphi_i$, which implies that
$\bar{\varphi}_i = \psi_i$ for $i=1,2$. By Lemma 1, $\bar{\varphi}
= \psi_1\otimes\psi_2 = \psi $ and condition (i) in the definition
of the extension property holds. If $A_1$ or $A_2$ is non-unital,
it is an immediate consequence of condition (ii)$'$ in the
definition of the extension property that condition (ii)$'$ holds
for $C_1\otimes C_2$ relative to $A_1\otimes_{\beta}A_2$. Thus
$C_1\otimes C_2$ has the extension property. Conversely, if
$C_1\otimes C_2$ has the extension property, it is a simple
exercise using similar methods to show that $C_2$ has the
extension property. \qed

\bigskip {\bf Note.} 1. An alternative proof of Theorem 2
when $A_1$ is unital can be given using the characterisation of
the extension property in terms of unitary conjugates given in the
paragraph following condition (ii)$'$ at the beginning of this
section.

2. If both $C_1$ and $C_2$ have the extension property a more
direct proof of the first part of Corollary 3 can be given as
follows. If $C$ is an abelian C*-subalgebra of
$A_1\otimes_{\beta}A_2$ containing $C_1\otimes C_2$, let
$\varphi=\varphi_1\otimes \varphi_2$ be a character of $C_1\otimes
C_2$. If $\bar{\varphi}$ is a character of $C$ extending
$\varphi$, then $\bar{\varphi}$ extends to a pure state $\psi$ of
$A_1\otimes_{\beta}A_2$ which equals $\psi_1\otimes \psi_2$, where
$\psi_i$ is the unique state extension of $\varphi_i$ to $A_i$ for
$i=1,2$, by the argument of the second paragraph of the above
proof. Thus $\bar{\varphi} = (\psi_1\otimes\psi_2)|_C$, which
means that $\varphi$ has a unique character extension to $C$.
Moreover no pure state of $A_1\otimes_{\beta}A_2$ and hence no
character of $C$ has a restriction to $C_1\otimes C_2$ equal to
$0$, since $C_1\otimes C_2$ contains an approximate identity for
$A_1\otimes_{\beta}A_2$. It follows by the Stone-Weierstrass
theorem that $C= C_1\otimes C_2$.

\bigskip {\bf 3. Masas with tensor products which are not maximal
abelian.}

\bigskip\noindent
{\bf 3.1.} Let $A=C^*_r({\Bbb F}_2) + K(\ell^2({\Bbb F}_2))$ in
$B(\ell^2({\Bbb F}_2))$, where $K(\ell^2({\Bbb F}_2))$ denotes the
compact linear operators on $\ell^2({\Bbb F}_2)$. Then $A$ is a
C*-algebra, $K(\ell^2({\Bbb F}_2))$ is an ideal of $A$ and
$A/K(\ell^2({\Bbb F}_2))\cong C^*_r({\Bbb F}_2)$. Let $q$ be the
canonical quotient map from $A$ onto $C^*_r({\Bbb F}_2)$, and let
$\lambda$ and $\rho$ denote the representations of $C^*_r({\Bbb
F}_2)$ corresponding to the left- and right-regular
representations of ${\Bbb F}_2$ on $\ell^2({\Bbb F}_2)$,
respectively. Then $\{\lambda\circ q, \rho\circ q\}$ is a
commuting pair of representations of the pair $\{A,A\}$ on
$\ell^2({\Bbb F}_2)$ with corresponding representation $\pi_r$ of
$A\odot A$ given by
   \[\pi_r(\sum a_i\otimes b_i) =
   \sum\lambda(q(a_i))\rho(q(b_i)).\]
A C*-norm $\|\phantom{x}\|_{\alpha}$ on $A\odot A$ is defined by
   \[\|x\|_{\alpha} = \max \{\|x\|_{min}, \|\pi_r(x)\|\}
   \quad(x\in A\odot A).\]
Let $\{\xi_g:g\in {\Bbb F}_2\}$ be the canonical orthonormal basis
of $\ell^2({\Bbb F}_2)$, for each $g\in {\Bbb F}_2$ let $e_g$ be
the projection onto the one dimensional subspace ${\Bbb C}\xi_g$
and let $C$ be the abelian C*-algebra generated by $\{e_g:g\in
{\Bbb F}_2\}\cup\{1\}$. Then $C\subset K(\ell^2({\Bbb F}_2))+
{\Bbb C}.1 \subset A$.

\bigskip
{\sc Proposition 4.} {\it The algebra $C$ is maximal abelian in
$A$, but $C\otimes C$ is not maximal abelian in
$A\otimes_{\beta}A$ for any C*-norm $\|\phantom{x}\|_{\beta}$
satisfying $\|x\|_\beta \geq \|x\|_\alpha$ on $x\in A\odot A$, in
particular if $\|\phantom{x}\|_\beta = \|\phantom{x}\|_{max}$.}

\bigskip
{\it Proof}. 1. To see that $C$ is maximal abelian in $A$, let $L$
be the closure of $C$ in the weak operator topology. Then $L$ is a
maximal abelian $*$-subalgebra of $B(\ell^2({\Bbb F}_2))$
isomorphic to $\ell^\infty({\Bbb N})$ and the canonical projection
$\sigma$ from $B(\ell^2({\Bbb F}_2))$ onto $L$ is given by
   \[\sigma(x) = \sum_{g\in {\Bbb F}_2}e_gxe_g\]
for $x \in B(\ell^2({\Bbb F}_2))$, the convergence of the sum on
the right being in the strong operator topology. If $x\in C^c =
C'\cap A$, then $\sigma(x) = \sum_{g\in {\Bbb F}_2}e_gxe_g =
\sum_{g\in {\Bbb F}_2}xe_g = x$. Now $x= k+a$, where $k\in
K(\ell^2({\Bbb F}_2))$ and $a\in C^*_r({\Bbb F}_2)$. Since
$e_g\lambda_h e_g = 0$ for $g, h \in {\Bbb F}_2$ when $h$ is not
the identity of ${\Bbb F}_2$, $\sigma(a) \in {\Bbb C}.1$ and since
$\sigma(k)\in C$, it follows that $\sigma(x)\in C$. Thus $C^c =
C$.

\bigskip
2. By [{\bf 11}] (see also [{\bf 13}]) the representation of
$C^*_r({\Bbb F}_2)\odot C^*_r({\Bbb F}_2)$ on $\ell^2({\Bbb F}_2)$
given by
   \[\sum a_i\otimes b_i \to \sum \lambda(a_i)\rho(b_i)\]
is not continuous relative to $\|\phantom{x}\|_{min}$, which
implies that there is a non-zero element in the kernel of the
canonical homomorphism from $A\otimes_{\alpha}A$ to
$A\otimes_{min}A$ and hence a non-zero element $x$ in the kernel
of the canonical homomorphism $\Psi$ from $A\otimes_{\beta}A$ to
$A\otimes_{min}A$. Since $K(\ell^2({\Bbb F}_2))$ is nuclear, the
restriction of $\|\phantom{x}\|_{\beta}$ to $K(\ell^2({\Bbb
F}_2))\odot A$ coincides with $\|\phantom{x}\|_{min}$ and
$\Psi|_{K(\ell^2({\Bbb F}_2))\otimes_{\beta}A}$ is isometric. Thus
for $k\in K(\ell^2({\Bbb F}_2))$,
   \[\|x(k\otimes 1)\|_{\beta} = \|\Psi(x(k\otimes 1))\|_{min} =
   \|\Psi(x)\Psi(k\otimes 1)\|_{min} = 0,\]
which implies that $x(k\otimes 1) = 0$. Similarly $(k\otimes 1)x =
0$. For $c\in C$ with $c= k + \lambda 1$, where $k\in
K(\ell^2({\Bbb F}_2))$ and $\lambda\in {\Bbb C}$,
   \[[x,c\otimes 1] = [x,k\otimes 1] = 0,\]
so that $x\in (C\otimes 1)^c$. Similarly $x\in (1\otimes C)^c$, so
that $x\in (C\otimes C)^c$. Since $\Psi(x)= 0$ and
$\Psi|_{C\otimes C}$ is isometric, $x\not\in C\otimes C$, which
means that $C\otimes C$ is not maximal abelian in
$A\otimes_{\beta}A$. \qed

\bigskip\noindent {\bf 3.2.} This result shows that without the requirement in the statement of
Corollary 3 that $C_1$ have the extension property, $C_1\otimes
C_2$ may not be maximal abelian in $A_1\otimes_{\beta}A_2$. By
modifying the construction of the C*-algebra $A$ above, it also
follows that the conclusion of the Corollary can fail if $C_2$
does not contain an approximate identity for $A_2$.

To see this, let $H$ be a separable infinite dimensional Hilbert
space with orthonormal basis $\{\xi_i:i\in {\Bbb N}\}$, and let
$H_1$ and $H_2$ be the closures of the linear subspaces of $H$
spanned by $\{\xi_{2i}:i\in {\Bbb N}\}$ and $\{\xi_{2i-1}:i\in
{\Bbb N}\}$, respectively. Then $H= H_1\oplus H_2$. A self-adjoint
unitary operator $u$ on $H$ is defined by
   \[u\xi_{2i} = {1\over \sqrt{2}}(\xi_{2i} + \xi_{2i-1}), \quad
   u\xi_{2i-1} = {1\over \sqrt{2}}(\xi_{2i} - \xi_{2i - 1}).\]
Let $\{e_{ij}\}$ be the set of rank one matrix units associated
with the basis $\{\xi_i\}$. Then
   \[u = \sum_i {1\over \sqrt{2}}(e_{2i,2i} + e_{2i, 2i-1} +
   e_{2i-1, 2i} - e_{2i-1, 2i-1}),\]
the sum on the right hand side converging in the strong operator
topology. If $t\in B(H)$ is such that $tH_1 \subset H_1$ and
$t|_{H_2} = 0$, then
   \[t = \sum_{i,j}t_{ij}e_{2i,2j},\]
where $t_{ij} = (t\xi_{2j}|\xi_{2i})$, the convergence again being
in the strong operator topology. Then
   \[utu = {1\over 2}\sum_{i,j}t_{ij}(e_{2i,2j} + e_{2i, 2j-1} + e_{2i-1, 2j} + e_{2i-1,
   2j-1}).\]
Let $B=C^*_r({\Bbb F}_2)\oplus \{0\}$ on $H$, where $C^*_r({\Bbb
F}_2)$ is acting in its identity representation on $H_1$ with
$H_1$ (respectively $H_2$) and $\ell^2({\Bbb F}_2)$ identified so
that $\{\xi_{2i}:i\in {\Bbb N}\}$ (respectively $\{\xi_{2i-1}:i\in
{\Bbb N}\}$) is the standard basis of $\ell^2({\Bbb F}_2)$ in some
enumeration, and let $A_0$ be the non-unital C*-algebra
$uBu+K(H)$. Then $A_0/K(\ell^2({\Bbb F}_2)) \cong C^*_r({\Bbb
F}_2)$.

For $t=\lambda_g\in C^*_r({\Bbb F}_2)$ and $i\in {\Bbb N}$,
   \[t_{ii}=(t\xi_{2i}|\xi_{2i})= \left\{\begin{array}{ll}
   1 & (g = e)\\
   0 & (g\neq e)
   \end{array}\right.,\]
where $e$ denotes the identity element of ${\Bbb F}_2$. Moreover
$e_{2i}u(t\oplus 0)ue_{2i} = {1\over 2}t_{ii}e_{2i,2i}$ and
$e_{2i-1}u(t\oplus 0)ue_{2i-1} = {1\over 2}t_{ii}e_{2i-1,2i-1}$
for $t\in C^*_r({\Bbb F}_2)$ and $i\in {\Bbb N}$. Let $C_0$ be the
abelian C*-algebra generated by the set of projections
$\{e_{ii}:i\in {\Bbb N}\}$ in $A_0$. The weak-operator closure $L$
of $C_0$ is maximal abelian in $B(H)$ and if $\sigma$ is the
canonical projection onto $L$, then for $t\in C^*_r({\Bbb F}_2)$,
   \[\sigma(u(t\oplus 0)u) = \sum_ie_iu(t\oplus 0)ue_i = {1\over
   2}\sum_it_{ii}(e_{2i,2i} + e_{2i-1,2i-1}) = \lambda 1\]
for some $\lambda\in {\Bbb C}$, where $1$ is the identity operator
on $H$. If $k+b\in C_0^c$ with $k\in K(H)$ and $b\in uBu$, then
$k+b\in L^c = L$ and so $k+b = \sigma(k) + \sigma(b)$. Now
$\sigma(k)\in C_0$ and $\sigma(b)=\lambda 1$ for some $\lambda\in
{\Bbb C}$. Hence $\lambda 1\in A_0$. Since $A_0$ is non-unital,
$\lambda=0$, which implies that $k+b=\sigma(k) \in C_0$. Thus
$C_0^c = C_0$ and $C_0$ is maximal abelian in $A_0$.

Let $q:A_0\to C^*_r({\Bbb F}_2)$ be the quotient map. With
$\lambda$ and $\rho$ as in 3.1, $\{\lambda, \rho\circ q\}$ is a
commuting pair of representations of the pair $\{C^*_r({\Bbb
F}_2), A_0\}$ on $\ell^2({\Bbb F}_2)$ with corresponding
representation $\pi_r$ of $C^*_r({\Bbb F}_2)\odot A_0$ given by
   \[\pi_r(\sum a_i\otimes b_i) =
   \sum\lambda(a_i)\rho(q(b_i)).\]
A C*-norm $\|\phantom{x}\|_{\alpha}$ on $C^*_r({\Bbb F}_2)\odot
A_0$ is defined by
   \[\|x\|_{\alpha} = \max \{\|x\|_{min}, \|\pi_r(x)\|\}
   \quad(x\in C^*_r({\Bbb F}_2)\odot A_0).\]
As in the proof of Proposition 4, for any norm
$\|\phantom{x}\|_{\beta}$ on $C^*_r({\Bbb F}_2)\odot A_0$ with
$\|x\|_{\beta}\geq \|x\|_{\alpha}$, there is a non-zero $x\in
C^*_r({\Bbb F}_2)\otimes_{\beta} A_0$ such that $\Psi(x) = 0$,
where $\Psi$ is the canonical homomorphism from  $C^*_r({\Bbb
F}_2)\otimes_{\beta} A_0$ to $C^*_r({\Bbb F}_2)\otimes_{min} A_0$.
Since $K(H)$ is nuclear, $\Psi|_{C^*_r({\Bbb
F}_2)\otimes_{\beta}K(H)}$ is isometric and, as in the proof of
Proposition 4, $[x,a\otimes k]=0$ for any $a\in C^*_r({\Bbb F}_2)$
and $k\in K(H)$. Thus $x\in (C_1\otimes C_0)^c\setminus(C_1\otimes
C_0)$ for any masa $C_1$ in $C^*_r({\Bbb F}_2)$, in particular if
$C_1$ is the abelian C*-subalgebra of $C^*_r({\Bbb F}_2)$ with the
extension property generated by one of the canonical unitary
generators. It is easy to see directly that $C_0$ does not contain
an approximate identity for $A_0$.

\bigskip {\bf 4. Some open problems.} When $H$ is the Hilbert space $\ell^2({\Bbb N})$, $B(H)$ is
non-nuclear and in particular the C*-norms $\|\phantom{x}\|_{max}$
and $\|\phantom{x}\|_{min}$ on $B(H)\odot C^*_r({\Bbb F}_2)$ are
distinct [{\bf 13}]. If $C$ is a non-atomic masa of $B(H)$
isomorphic to $L^\infty(0,1)$ (which does not have the extension
property by [{\bf 9}]), is it true that $(C\otimes 1)^c = C\otimes
C^*_r({\Bbb F}_2)$ in $B(H)\otimes_{max}C^*_r({\Bbb F}_2)$? Junge
and Pisier [{\bf 7}] have shown that the C*-norms
$\|\phantom{x}\|_{max}$ and $\|\phantom{x}\|_{min}$ on $B(H)\odot
B(H)$ are distinct. Is it true that $(C\otimes 1)^c = C\otimes
B(H)$ in $B(H)\otimes_{max}B(H)$ for any masa $C$ of $B(H)$? Is
$C_1\otimes C_2$ maximal abelian in $B(H)\otimes_{max} B(H)$ for
any masas $C_1$ and $C_2$? The case $C_1\cong C_2\cong
\ell^\infty({\Bbb N})$ is particularly intriguing. The question of
whether masas isomorphic to $\ell^\infty({\Bbb N})$ have the
extension property relative to $B(H)$ was first investigated by
Kadison and Singer [{\bf 9}], but remains at the time of writing
one of the more significant unsolved problems in the subject,
despite the attention of many distinguished workers. If the
Kadison-Singer problem had a positive solution, it would follow by
Corollary 3 that $\ell^\infty({\Bbb N})\otimes \ell^\infty({\Bbb
N})$ is maximal abelian in $B(H)\otimes_{max}B(H)$. If however it
could be shown that $\ell^\infty({\Bbb N})\otimes
\ell^\infty({\Bbb N})$ is not maximal abelian in
$B(H)\otimes_{max}B(H)$, then it would follow that
$\ell^\infty({\Bbb N})$ does not have the extension property.

\bigskip
{\sc Acknowledgement}. The author is grateful to Stuart White for
a stimulating discussion on the topics considered here and to the
referee for pointing out a number of inaccuracies.

\bigskip
\centerline{REFERENCES}

\bigskip
\noindent \small
\parindent=.55in

{\bf 1.} C.A.\ Akemann, and F.\ Shultz, Perfect $C\sp *$-algebras,
{\it Mem.\ Amer.\ Math.\ Soc.} {\bf 55} (1985), no.\ 326

{\bf 2.} J.\ Anderson, A maximal abelian subalgebra of the Calkin
algebra with the extension property, {\it Math.\ Scand.} {\bf 42}
(1978), 101--110.

{\bf 3.} R.J.\ Archbold, On the centre of a tensor product of
$C\sp*$-algebras, {\it J.\ London Math.\ Soc.} (2) {\bf 10}
(1975), 257--262.

{\bf 4.} R.J.\ Archbold, Extensions of states of $C\sp{*}
$-algebras, {\it J.\ London Math.\ Soc.} (2) {\bf 21} (1980),
351--354.

{\bf 5.} C.J.K.\ Batty, On relative commutants in tensor products
of $C\sp{\sp*}$-algebras, {\it Math.\ Z.} {\bf 151} (1976),
215--218.

{\bf 6.} R.\ Haydon and S.\ Wassermann, A commutation result for
tensor products of $C\sp{*} $-algebras, {\it Bull.\ London Math.\
Soc.} {\bf 5 }(1973), 283--287.

{\bf 7.} M.\ Junge and G.\ Pisier, Bilinear forms on exact
operator spaces and $B(H)\otimes B(H)$, {\it Geom.\ Funct.\ Anal.}
{\bf 5} (1995), 329--363.

{\bf 8.} C.\ Lance, Tensor products of non-unital
$C\sp*$-algebras, {J.\ London Math.\ Soc.} (2) {\bf 12} (1975/76),
160--168.

{\bf 9.} R.\ Kadison and I.M.\ Singer, Extensions of pure states,
{\it Amer.\ J.\ Math.} {\bf 81} (1959), 383--400.

{\bf 10.} S.\ Sakai, {\it $C\sp*$-algebras and $W\sp*$-algebras},
Ergebnisse der Mathematik und ihrer Grenzgebiete, Band 60,
Springer, 1971.

{\bf 11.} M.\ Takesaki, On the cross-norm of the direct product of
$C\sp{*} $-algebras, {\it T\^ohoku Math.\ J.} (2) {\bf 16} (1964),
111--122.

{\bf 12.} S.\ Wassermann, The slice map problem for
$C\sp*$-algebras, {\it Proc.\ London Math.\ Soc.} (3) 32 (1976),
537--559.

{\bf 13.} S.\ Wassermann, On tensor products of certain group
$C\sp{*} $-algebras, {\it J.\ Functional Analysis} {\bf 23}
(1976), 239--254.

{\bf 14.} S.\ Wassermann, Tensor products of free-group
C*-algebras, {\it Bull.\ London Math. Soc.} {\bf 22} (1990),
375--380.

\end{document}